\documentclass[a4paper,leqno]{amsart}
\usepackage[utf-8]{inputenc}
\usepackage[english]{babel}
\usepackage{latexsym,amsthm,amsmath,amscd,amsxtra,amssymb,amsxtra,exscale,mathrsfs,pgf}
\usepackage[colorlinks,linkcolor={blue},citecolor={blue},urlcolor={red},]{hyperref}

\theoremstyle{plain}
\newtheorem{theorem}{Theorem}[section]
\newtheorem{corollary}[theorem]{Corollary}
\newtheorem{lemma}[theorem]{Lemma}
\newtheorem{proposition}[theorem]{Proposition}

\theoremstyle{definition}
\newtheorem{definition}[theorem]{Definition}

\theoremstyle{remark}

\newtheorem*{remark*}{Remark}

\colorlet{shaded}{black!25!white}

\newcommand{\RR}{{\mathbb R}}
\newcommand{\CC}{{\mathbb C}}

\newcommand{\al}{\alpha}
\newcommand{\ga}{\gamma}

\newcommand{\eps}{\varepsilon}
\newcommand{\la}{\lambda}
\newcommand{\om}{\omega}

\newcommand{\si}{\sigma}
\newcommand{\Laplace}{\mathcal L}

\newcommand{\cD}{{\mathcal D}}

\newcommand{\norm}[1] {\| #1 \|}

\newcommand{\bignorm}[1]{\bigl\| #1 \bigr\|}
\newcommand{\Bignorm}[1]{\Bigl\| #1 \Bigr\|}

\newcommand{\sfrac}[2] { {{}^{#1}\!\!/\!{}_{#2}}}
\newcommand{\einhalb} {\sfrac{1}{2}}
\newcommand{\pihalbe} {\sfrac{\pi\,}{2}}

\renewcommand{\Re}[1] {{\rm Re}\,#1}

\newcommand{\suchthat}{\,\colon\;}
\newcommand{\Sec}[1]{\Sigma_{#1}}

\newcommand {\comment}[1]{\relax}

\newcommand{\loc}{\text{loc}}

\newcommand{\ahut}{\widehat{a}}
\newcommand{\logplus}{\log^{+}\!}
\newcommand{\REGULAR}[1]{C^{\text{reg}}_{#1}}
\newcommand{\BOUNDED}{{\mathcal B}}

\begin{document}

  \title{Control of Volterra systems with scalar kernels}

  \author[Bernhard H. Haak]{Bernhard H. Haak}
  \address{Institut de Math\'ematiques de Bordeaux\\
           351 cours de la  Lib\'eration\\33405 Talence CEDEX\\France}
  \email{Bernhard.Haak@math.u-bordeaux1.fr}

  \author[Birgit Jacob]{Birgit Jacob}
  \address{Institut f\"ur Mathematik\\ Universit\"at Paderborn\\ Warburger Stra\ss e 100\\D-33098 Paderborn\\Germany}
  \email{jacob@math.uni-paderborn.de}

  \date\today

  \begin{abstract}
     \noindent Volterra observations systems with scalar kernels are
     studied. New sufficient conditions for admissibility of observation
     operators are developed. The obtained results are applied to
     time-fractional diffusion equations of distributed order.
   \end{abstract}

  \maketitle

\allowdisplaybreaks

\section{Introduction}

Consider the following scalar abstract Volterra system
\begin{equation} \label{eq:volterra-uncontrolled}
    x(t) = f(t) + \displaystyle\int_0^t a(t{-}s)\, A x(s)\, ds.
\end{equation}
Here, the operator $A$ is supposed to be a closed operator with dense
domain on a Banach space $X$ having its spectrum contained in some
open sectorial region of the complex plane, symmetric to the real axis
and open to the left:
\[
\si(A) \subseteq -\Sec{\om}
       \qquad\text{where}\qquad \Sec{\om}  = \left\{ z\in\CC
         \suchthat |\arg(z)| < \om \right\}
\]
for some $\om \in (0,\pi)$. Moreover, the resolvent of $A$ is supposed
to satisfy a growth condition of the type $\norm{\la R(\la, A) } \le
M$ uniformly on each sector $\Sec{\pi{-}\om{-}\eps}$.  Typical
examples of such operators are generators of bounded strongly
continuous semigroups, where $\om \le \pihalbe$. We call $-A$ a
{\em  sectorial operator of type $\om \in (0, \pi)$}, but we mention
that 'sectoriality' may have different meanings for different authors
in the literature.

\smallskip

The kernel function $a \in L^1_{\text{loc}}$ is supposed to be of
sub-exponential growth so that its Laplace transform $\ahut(\la)$
exists for all $\la$ with positive real part. The kernel is called
{\em sectorial of angle $\theta \in (0,\pi)$} if
\[
\ahut(\la) \in \Sec{\theta} \qquad \text{for all }\la\text{ with
  positive real part.}
\]
We will consider only {\em parabolic} equations
(\ref{eq:volterra-uncontrolled}) in the sense of Pruess
\cite{Pruess:Volterra}. In the case that $-A$ is sectorial of some
angle $\om\in(0,\pi)$ this is equivalent to require $\ahut(\la)
\not=0$ and $\tfrac1{\ahut(\la)} \in \varrho(A)$ for all $\la$ with
positive real part.

In particular, when $-A$ and $a$ are both sectorial in the
respective sense with angles that sum up to a constant strictly
inferior to $\pi$, the Volterra equation is parabolic.

\smallskip

The kernel function is said to be {\em $k$-regular} if there is a constant
$K>0$ such that
\[
       |\la^n \ahut^{(n)}(\la) | \le K |\ahut(\la)|
\]
for all $n=0,1\ldots k$ and all $\la$ with positive real part. In Pruess
\cite[Theorem I.3.1]{Pruess:Volterra} it is shown that parabolic
equations with a $k$-regular kernel for $k\ge 1$ admit a unique {\em solution
family}, i.e.~a  family of bounded linear operators $(S(t))_{t\ge 0}$
on $X$, such that
\begin{enumerate}
\item $S(0)=I$ and  $S(\cdot)$ is strongly continuous on $\mathbb R_+$.
\item $S(t)$ commutes with $A$, which means $S(t)(D(A))\subset D(A)$ for all
  $t\ge 0$, and $AS(t)x=S(t)Ax$ for all $x\in D(A)$  and $t\ge 0$.
\item For all $x\in D(A)$ and all $t\ge 0$ the resolvent equations hold:
\begin{equation}\label{eqn:resolvent}
S(t)x = x +\int_0^t a(t-s)AS(s)x\,ds.
\end{equation}
\end{enumerate}
Moreover, $S \in C^{k-1}((0,\infty), \BOUNDED(X) )$ and $\norm{ t^n
  S^{(n)}(t) } \le K$ for all $n=0,\ldots, k-1$.

\bigskip

The purpose of this article is to present conditions for the
admissibility of observation operators to parabolic Volterra
equations, that is, we consider the 'observed' system
\begin{equation} \label{eq:volterra-system} \tag{V}
\left\{
  \begin{split}
     x(t) &= f(t) + \int_0^t a(t{-}s)\, A x(s)\, ds \\
     y(t) &= C x(t)
  \end{split}
\right.
\end{equation}

Additionally to the sectoriality condition on $A$ and the
parabolicity condition on the Volterra equation, the operator $C$ in
the second line is supposed to be an operator from $X$ into another
Banach space $Y$ that acts as a bounded operator from $X_1 \to Y$
where $X_1 = \cD(A)$ endowed by the graph norm of $A$. In order to
guarantee that the output function lies locally in $L_2$ we impose
the following condition.

\begin{definition}\label{def:volterra-admiss}
  A bounded linear operator $C: X_1 \to Y$ is called {\em finite-time
    admissible} for the Volterra equation
  (\ref{eq:volterra-uncontrolled}) if there are constants $\eta, K>0$
  such that
\[
   \biggl(\int_0^t \norm{ C S(r) x}^2 \,dr\biggr)^\einhalb
\le K e^{\eta t} \norm{x}
\]
for all $t\ge 0$ and all $x\in \cD(A)$.
\end{definition}

The notion of admissible observation operators is well studied in
the literature for Cauchy systems, that is, $a\equiv 1$, see for example
\cite{JacobPartington:survey}, \cite{Staffans}, and
\cite{Weiss:admiss-observation}. Admissible observation operators
for Volterra systems are studied in \cite{hjpp08},
\cite{JacobPartington04volterra}, \cite{japa07} and \cite{Jung}.

The Laplace transform of $S$, denoted by $H$, is given by
\[ H(\lambda)x=\frac{1}{\lambda}(I-\hat{a}(\lambda)A)^{-1} x, \qquad
{\rm Re}\,\lambda>0.\]

The following necessary condition for admissibility was shown in
\cite{japa07}.

\begin{proposition}\label{lem:bdd-condition}
If $C$ is a finite-time admissible observation operator for the
  Volterra equation (\ref{eq:volterra-uncontrolled}), then there is a
  constant $M>0$ such that
\begin{equation}\label{eqnbdd}
    \norm{ \sqrt{{\rm Re}\,\lambda} C H(\lambda)}
\le     M , \qquad {\rm Re}\,\lambda>0.
\end{equation}
\end{proposition}

In \cite{japa07} it is shown that (\ref{eqnbdd}) is also sufficient
for admissibility if $X$ is a Hilbert space, $Y$ is finite-dimensional
and $A$ generates a contraction semigroup.  However, in general this
condition is not sufficient (see e.g. \cite{JacobPartington:survey}).

We show that the slightly stronger growth condition on the resolvent
\begin{equation*} 
   \sup_{r>0}  \bignorm{ \,\ (1+\logplus r)^\al r^\einhalb C H(r) \,}
  < \infty,
\end{equation*}
is sufficient for admissibility if $\al>\einhalb$ (see
Theorem~\ref{thm:hans}). This result generalizes the sufficient
condition of Zwart \cite{Zwart:logplus} for Cauchy systems to general Volterra
systems (\ref{eq:volterra-uncontrolled}).

Our second main result, Theorem~\ref{thm:theo1} provides a perturbation
argument to obtain admissibility for the controlled Volterra equation
from the admissibility of the control operator for the underlying
Cauchy equation. In the particular case of diagonal semigroups and
one-dimensional output spaces $Y$ this improves a direct Carleson
measure criterion from Haak, Jacob, Partington and Pott \cite{hjpp08}.

We proceed as follows. In Section 2 we obtain an integral
representation for the solution family $(S(t))_{t\ge 0}$ and several
regularity results of the corresponding kernel. Section 3 is devoted
to sufficient condition for admissibility of observation operators.  A
perturbation result as well as a general sufficient condition is
obtain. Several examples are included as well.

To enhance readability of the calculations, for rest of this article,
$K$ denotes some positive constant that may change from one line to
the other unless explicitly quantified.

\section{Regularity transfer}

The main result of this section is formulated in the following
proposition. Let $s(t, \mu)$ denote the solution of the scalar
equation
\[
    s(t, \mu) + \mu \int_0^t a(t{-}r) s(r, \mu)\,dr = 1 \qquad t>0, \mu\in \CC.
\]

\begin{proposition}\label{prop2.3}
  Suppose that the kernel $a \in L^1_\loc(\RR_+)$ is $1$-regular, and
  sectorial of angle $\theta < \pihalbe$. Then there exists a family
  of functions $v_t$ such that
\[
\Laplace(v_t)(\mu) = s(t, \mu) \quad \text{and}\quad
S(t) = \int_0^\infty v_t(s) T(s)\,ds
\]
satisfying
\begin{enumerate}
\item $\sup_{t>0} \norm{ v_t }_{L^1(\RR_+)} < \infty$
\item $\norm{ v_t }_{L^2(\RR_+)} \le K (t^{-\sfrac\theta\pi} +
  t^{+\sfrac\theta\pi})$ where
    $K$ depends only on $\theta$ and $\REGULAR{a, 1}$.
\item $\norm{ v_t }_{W^{1,1}} \le K (1+t^{-\frac{2\theta}\pi} + t^{+\frac{2\theta}\pi})$.
\end{enumerate}
Moreover, for $\ga\in [0,1]$, $|\mu^\ga s(t, \mu)| \le K t^{-\frac{2\ga\theta}\pi}$
\end{proposition}

For the proof of this proposition the following two lemmas are
needed.

\begin{lemma} \label{lem:automatisches-wachstum}
  Suppose $a\in L^1_{loc}(\mathbb R_+)$ is $1$-regular and sectorial
  of angle $\theta\le\pi$.  Let $\rho_0:=2\theta/\pi$.  Then there
  exists a constant $c>0$ such that
\[
   |\ahut(\lambda)| \ge
  \begin{cases}
    c|\lambda|^{-\rho_0} &|\lambda|\ge 1\\
    c|\lambda|^{\rho_0} &|\lambda| \le 1
  \end{cases}
\]
for all $\lambda \in \CC_+$.
\end{lemma}
\begin{proof}
  We borrow the argument from the proof of \cite[Proposition
  1]{MonniauxPruess}: we start with the analytic completion of the
  Poisson formula for the harmonic function
  $H(\lambda)=\arg\ahut(\lambda)$, that is,
\[
  \log \ahut(\lambda) =\kappa_0 +\frac{i}{\pi}\int_{-\infty}^\infty
  \left[\frac{1-i\rho\lambda}{\lambda-i\rho}\right]h(i\rho)\frac{d\rho}{1+\rho^2},
\]
where $\kappa_0\in\mathbb R$ is a constant.  An easy calculation shows
\[
  |\Re \log\ahut(\lambda)|\le \kappa_0+\rho_0 |\log \lambda|
\]
for real $\lambda>0$, and thus
\[
  |\ahut(\lambda)| = e^{\log(|\ahut(\la)|)} = e^{\Re \log \ahut(\la)} \ge
   \begin{cases} c\lambda^{-\rho_0} &\lambda\ge 1\\
                  c\lambda^{\rho_0} & 0\le \lambda \le 1
   \end{cases},
\]
where $c:=e^{-\kappa_0}>0$. This estimate, together with \cite[Lemma
8.1]{Pruess:Volterra} stating the existence of a constant $c>0$ such
that $c^{-1} \leq \bigl|\ahut(|\la|) / \ahut(\la) \bigr| \leq c$ for
all $\la \in \CC_+$ completes the proof.
\end{proof}

\begin{lemma}\label{lem:umgekehrte-dreicksungleichung}
  Let $\theta \in (0, \pi)$. Then there exists $c_\theta >0$ such that
  \begin{equation}   \label{eq:umgekehrte-dreicksungleichung}
      1+|\la| \le c_\theta \; |1+\la|
  \end{equation}
  for all $\la \in \Sec{\pi{-}\theta}$.
\end{lemma}
\begin{proof}
   \begin{figure}[h]\label{bild1}
     \centering
    \begin{pgfpicture}{0cm}{0.5cm}{7cm}{8.0cm}
    \pgfsetlinewidth{0.8pt}
    \color{white}
    \pgfmoveto{\pgfxy(-1,1)}
   \pgflineto{\pgfxy(-1,7)}
    \pgflineto{\pgfxy(8,7)}
    \pgflineto{\pgfxy(8,1)}
    \pgfclosepath
   \pgffill
    \color{shaded}
    \pgfmoveto{\pgfxy(3,3)}
    \pgflineto{\pgfxy(-1.1,4.4)}
    \pgflineto{\pgfxy(-1.1,1.6)}
    \pgfclosepath
    \pgffill
    \color{black}
    \pgfcircle[stroke]{\pgfxy(6.35,3)}{1pt}
    \pgfcircle[stroke]{\pgfxy(0.9,5.6)}{1pt}
    \pgfcircle[stroke]{\pgfxy(-0.15,4.12)}{1pt}
    \pgfcircle[stroke]{\pgfxy(1.8,3)}{1pt}
   \pgfsetdash{{.1cm}}{0pt}
    \pgfxyline(6.35,3)(0.9,5.6)
    \pgfxyline(3,3)(0.9,5.6)
    \pgfsetdash{}{0pt}
    \pgfmoveto{\pgfxy(1.5,3)}
    \pgfcurveto{\pgfxy(1.5,3.3)}{\pgfxy(1.55,3.3)}{\pgfxy(1.6,3.5)}
    \pgfmoveto{\pgfxy(4.2,3)}
    \pgfcurveto{\pgfxy(4.2,3.4)}{\pgfxy(4.25,3.7)}{\pgfxy(4.35,3.95)}
    \pgfmoveto{\pgfxy(3.7,3)}
    \pgfcurveto{\pgfxy(3.7,3.3)}{\pgfxy(3.73,3.25)}{\pgfxy(3.77,3.44)}
    \pgfmoveto{\pgfxy(-0.15,4.12)}
    \pgfcurveto{\pgfxy(0.9,7.2)}{\pgfxy(6,7.3)}{\pgfxy(6.35,3)}
    \pgfxyline(6.35,3)(-0.15,4.12)
    \pgfxyline(3,3)(-1,4.4)
    \pgfxyline(3,3)(-1,1.6)
    \pgfsetendarrow{\pgfarrowto}
    \pgfxyline(-1,3)(8,3)
    \pgfxyline(3,1)(3,7)
    \pgfputat{\pgfxy(6.35,2.65)}{\pgfbox[center,base]{$|\lambda|$}}
    \pgfputat{\pgfxy(0.8,5.7)}{\pgfbox[center,base]{$\lambda$}}
    \pgfputat{\pgfxy(-0.3,4.3)}{\pgfbox[center,base]{$\tilde\lambda$}}
    \pgfputat{\pgfxy(1.7,3.1)}{\pgfbox[center,base]{$\theta$}}
    \pgfputat{\pgfxy(1.8,2.7)}{\pgfbox[center,base]{$-1$}}
    \pgfputat{\pgfxy(4.6,3.5)}{\pgfbox[center,base]{$\alpha$}}
    \pgfputat{\pgfxy(4,3.05)}{\pgfbox[center,base]{$\tilde\alpha$}}
     \end{pgfpicture}
     \caption{Illustration of (\ref{eq:umgekehrte-dreicksungleichung})}
 \end{figure}
Clearly, $\al >\widetilde \al$, see Figure 1. Since $\widetilde \al
= \tfrac{\theta}2$, the assertion follows then from the fact that
$\frac{|1+\la|}{1+|\la|} = \frac{\sin(\al)}{\sin(\theta{-}\al)} \ge
\sin(\al) \ge \sin(\theta/2)$.
\end{proof}

\begin{proof}[Proof of Proposition~\ref{prop2.3}]
\begin{enumerate}
\item is \cite[Proposition I.3.5]{Pruess:Volterra}. This latter result
  is also the principal inspiration of the next part:
\item
  Let $\si(\la, \mu) = \bigl(\Laplace s(\cdot, \mu)\bigr) (\la)$, i.e.
  $\si(\la,\mu) = \frac1{\la(1+\mu \ahut(\la) )}$. Fix $t>0$ and
  $\eps>0$. Then
\[
   s(t, \mu) =  \tfrac1{2\pi i} \int_{\eps-i\infty}^{\eps+i\infty} e^{\la t} \si(\la, \mu)\,d\la.
\]
Then, by partial integration
\begin{eqnarray*}
  s(t, \mu)
&= & \lim_{R\to \infty} \tfrac{1}{2\pi i} \biggl[
        \tfrac1t e^{\la t} \si(\la, \mu)
        \biggr]_{\la=\eps-iR}^{\la=\eps+iR}
    - \tfrac{1}{2\pi i} \int_{\eps-iR}^{\eps+iR} \tfrac1t e^{\la t}
        \frac{d}{d\la} \si(\la, \mu) \,d\la\\
&= & - \tfrac{1}{2\pi i} \int_{\eps-i\infty}^{\eps+i\infty} \tfrac1t e^{\la t}
        \frac{d}{d\la} \si(\la, \mu) \,d\la\\
\end{eqnarray*}
An elementary calculation gives
\[
   \frac{d}{d\la} \frac{1}{\la(1+\mu \ahut(\la))}
=  - \frac{1+\mu \ahut(\la)\Bigl(1+ \Bigl( \frac{\la \ahut'(\la)}{\ahut(\la)}\Bigr)\Bigr)}
      {\la^2 (1+\mu \ahut(\la))^2}
\]
By 1-sectoriality of the kernel, $\Bigl| \frac{\la
  \ahut'(\la)}{\ahut(\la)}\Bigr| \le C = \REGULAR{a, 1}$ and so
Lemma~\ref{lem:umgekehrte-dreicksungleichung} yields for any $\delta>0$,
\begin{eqnarray*}
  &&   \biggl(\int_{-\infty}^\infty \Bigl| s(t, \delta+iy)\Bigr|^2 \,dy\biggr)^\einhalb \\
&\le& C_\theta (1+C) \frac{e^{\eps t}}{2\pi t} \biggl(\int_{-\infty}^\infty \biggl(\int_{-\infty}^\infty
      \frac{1}{(\eps^2{+}x^2)(1+|\delta{+}iy| |\ahut(\eps{+}ix)|)} \,dx \biggr)^2 \,dy\biggr)^\einhalb\\
&\le& \sqrt{2} C_\theta (1+C) \frac{e^{\eps t}}{2\pi t} \biggl(\int_{0}^\infty \biggl(\int_{-\infty}^\infty
      \frac{1}{(\eps^2{+}x^2)(1+|y| |\ahut(\eps{+}ix)|)} \,dx \biggr)^2 \,dy\biggr)^\einhalb\\
&\le& \sqrt{2} C_\theta (1+C) \frac{e^{\eps t}}{2\pi t} \int_{-\infty}^\infty   \biggl(\int_{0}^\infty
      \frac{1}{(\eps^2{+}x^2)^2(1+|y| |\ahut(\eps{+}ix)|)^2} \,dy \biggr)^\einhalb \,dx\\
& = & \sqrt{2} C_\theta (1+C) \frac{e^{\eps t}}{2\pi t} \int_{-\infty}^\infty \frac1{(\eps^2{+}x^2)|\ahut(\eps+ix)|^\einhalb}  \biggl(\int_{0}^\infty  \frac{1}{(1+u)^2} \,du \biggr)^\einhalb \,dx\\
& = & \sqrt{2} C_\theta (1+C) \frac{e^{\eps t}}{2\pi t} \int_{-\infty}^\infty \frac1{(\eps^2{+}x^2)|\ahut(\eps+ix)|^\einhalb}   \,dx\\
\end{eqnarray*}
Now we split the integral into two parts, by considering the cases
$\eps^2{+}x^2 \ge 1$ and $\eps^2{+}x^2 < 1$ to apply
Lemma~\ref{lem:automatisches-wachstum} which is controlling $|1/\ahut|$. 
Substituting $x=\eps t$ in both parts easily gives
\[
\norm{ s(t, \cdot) }_{H^2} \le \widetilde C_\theta \frac{e^{\eps t}}{
  t}  \bigl(\eps^{-1-\sfrac\theta\pi} + \eps^{-1+\sfrac\theta\pi}  \bigr),
\]
which yields the assertion by letting $\eps = \sfrac1t$.
\item
We argue in the same spirit as above: by partial integration
\begin{eqnarray*}
 \tfrac{d}{d\mu} \Bigl(\mu s(t, \mu)\Bigr)
&= &  \tfrac{1}{2\pi i} \int_{\eps-i\infty}^{\eps+i\infty} \tfrac1t e^{\la t}
        \frac{d^2}{d\mu\,d\la} \Bigl( \mu \si(\la, \mu) \Bigr)\,d\la\\
\end{eqnarray*}
An elementary calculation gives
\[
   \frac{d^2}{d\la\,d\mu} \frac{\mu \ahut(\la)}{(\la(1+\mu \ahut(\la)))^2}
=   \frac{1+\mu \ahut(\la)\Bigl(1+ 2 \Bigl( \frac{\la \ahut'(\la)}{\ahut(\la)}\Bigr)\Bigr)}
      {\la^2 (1+\mu \ahut(\la))^3}
\]
By 1-sectoriality of the kernel, $\Bigl| \frac{\la
  \ahut'(\la)}{\ahut(\la)}\Bigr| \le C$ and so the Lemma yields any
$\delta>0$,
\begin{eqnarray*}
  &&   \int_{-\infty}^\infty \Bigl| \frac{d}{d\mu} \Bigl(\mu s(t, \delta+iy)\Bigr)\Bigr| \,dy\\
&\le& C_\theta (1+2C) \frac{e^{\eps t}}{2\pi t} \int_{-\infty}^\infty \int_{-\infty}^\infty
      \frac{1}{(\eps^2{+}x^2)(1+|\delta{+}iy| |\ahut(\eps{+}ix)|)^2} \,dx\,dy\\
&\le& C_\theta (1+2C) \frac{e^{\eps t}}{2\pi t} \int_{-\infty}^\infty \int_{-\infty}^\infty
      \frac{1}{(\eps^2{+}x^2)(1+|y| |\ahut(\eps{+}ix)|)^2} \,dx\,dy\\
&=&   2 C_\theta (1+2C) \frac{e^{\eps t}}{2\pi t}
       \int_{-\infty}^\infty
      \frac{1}{(\eps^2{+}x^2)}\frac{1}{|\ahut(\eps{+}ix)|}  \int_{0}^\infty \frac{1}{(1+u)^2} \,du \,dx \\
&=&   C_\theta (1+2C) \frac{e^{\eps t}}{2\pi t}
       \int_{-\infty}^\infty
      \frac{1}{(\eps^2{+}x^2)}\frac{1}{|\ahut(\eps{+}ix)|} \,dx \\
&\le& K (t^{-\frac{2\theta}\pi} + t^{+\frac{2\theta}\pi})
\end{eqnarray*}
by choosing $\eps=\sfrac1t$. This shows that
$f_t(\mu) = \tfrac{d}{d\mu} \Bigl(\mu s(t, \mu)\Bigr)  \in H^1(\CC_+)$.
We may apply Hardy's inequality (see e.g. \cite[p.198]{Duren:Hp-spaces},
\cite[Theorem 4.2]{HilleTamarkin}),
\[
    \int_0^\infty \frac{| \check{f_t}(r) |}{r}\,dr
\le \tfrac12 \int_{-\infty}^{\infty} | f_t(i\om)|\,d\om
\]
so that $\frac{\check{f_t}(r)}{r} \in L^1(\RR_+)$ is Laplace transformable
for every $t>0$. Since
\[
  \Laplace\Bigl(\frac{\check{f_t}(r)}{r}\Bigr)(\si)
= \int_\si^\infty f_t(\mu)\,d\mu = \si s(t, \si),
\]
we find that $\mu \mapsto \mu s(t, \mu) \in H^\infty(\CC_+)$ with a
norm controlled by a multiple of $(t^{\frac{+2\theta}\pi} +
t^{-\frac{2\theta}\pi})$. This implies that $v_t' \in L^1(\RR_+)$.
Together with (a) the claim follows.
  \end{enumerate}

Finally, the same technique gives an estimate for the growth of $s(t, \mu)$:
  \begin{eqnarray*}
     \mu^\ga s(t, \mu)
&\le& K \frac{|\mu|^\ga e^{\eps t}}{t}
       \int_{-\infty}^\infty \frac{1}{(\eps^2+r^2)(1+|\mu||\ahut(\eps+ir)|)}\,dr \\
&\le& K \frac{e^{\eps t}}{t}
       \int_{-\infty}^\infty \frac{1}{(\eps^2+r^2)|\ahut(\eps+ir)|^\ga} \,
        \frac{|\mu|^\ga |\ahut(\eps+ir)|^\ga}{(1+|\mu||\ahut(\eps+ir)|)}\,dr \\
&\le& K \frac{e^{\eps t}}{t}
       \int_{-\infty}^\infty \frac{1}{(\eps^2+r^2)|\ahut(\eps+ir)|^\ga} \,dr \\
&\overset{\eps=\frac1t}{\le}&
      K (t^{-\frac{2\ga\theta}\pi} + t^{+\frac{2\ga\theta}\pi}).
  \end{eqnarray*}
\end{proof}

\section{Sufficient conditions for finite-time admissibility}

In this section we present the two main results of this paper.

\begin{theorem}\label{thm:theo1}
  Let $A$ generate an exponentially stable strongly continuous
  semigroup $(T(t))_{t\ge 0}$ and let $C: X_1 \to Y$ be bounded. Further
  we assume that the kernel $a\in L^1_{loc}(\mathbb R_+)$ is
  $1$-regular and sectorial of angle $\theta<\pi/2$.  Then finite-time
  admissibility of $C$ for the semigroup $(T(t))_{t\ge 0}$ implies that
  of $C$ for the solution family $(S(t))_{t\ge 0}$.
\end{theorem}
\begin{proof}
  By Proposition~\ref{prop2.3} there exists a family of functions
  $v_t$ such that $\norm{v_t}_{L^2(\mathbb R_+)}\le K t^{-\theta/\pi}$
  for some constant $K>0$ independent of $t\ge 0$ and
\[
S(t) = \int_0^\infty v_t(r)\,T(r) \,dr, \qquad t\ge 0.
\]
For  $x\in \cD(A)$ we have thus
\[
C S(t)x = \int_0^\infty v_t(r)\, C T(r)x \,dr.
\]
Note that finite-time admissibility of $C$ for $(T(t))_{t\ge 0}$ implies
the existence of a constant $M>0$ such that
\[
   \norm{C T(\cdot)x}_{L^2(\mathbb R_+)} \le M\norm{x}, \qquad x\in \cD(A),
\]
thanks to the exponential stability of $(T(t))_{t\ge 0}$. Thus the
result follows from Cauchy-Schwarz inequality.
\end{proof}

By replacing the Cauchy-Schwarz inequality by H\"older's inequality,
similar arguments can be used to obtain sufficient conditions for
$L^p$-admissibility.

\begin{corollary}
Assume in addition to the hypotheses of the theorem that
 one of the following conditions is satisfied:
\begin{enumerate}
\item $Y$ is finite-dimensional, $X$ is a Hilbert space and $A$
  generates a contraction semigroup;
\item $X$ is a Hilbert space and $A$ generates a normal, analytic
  semigroup;
\item $A$ generates an analytic semigroup and $(-A)^{1/2}$ is an
  finite-time admissible observation operator for $(T(t))_{t\ge 0}$.
\end{enumerate}
If there exists a constant $M>0$ such that
\begin{equation}\label{eqnresolvent}
  \norm{C(\lambda-A)^{-1}} \le \frac{M}{\sqrt{{\rm Re}\,\lambda}}, \qquad
  {\rm Re}\,\lambda>0,
\end{equation}
then $C$ is a finite-time admissible observation operator for
$(S(t))_{t\ge 0}$.
\end{corollary}
\begin{proof}
  Under the assumption of the corollary, the inequality
  (\ref{eqnresolvent}) implies that $C$ is a finite-time admissible
  observation operator for $(S(t))_{t\ge 0}$, see \cite{JacobPartington:contraction},
  \cite{HansenWeiss:Carleson}, \cite{LeMerdy:weiss-conj}.  Thus the result
  follows from Theorem~\ref{thm:theo1}.
\end{proof}

The following corollary is an immediate consequence of the Carleson
measure criterion of Ho and Russell \cite{HoRussell}.

\begin{corollary}\label{cor:carleson}
Assume in addition to the hypotheses of the theorem that
$A$ admits a Riesz basis of eigenfunctions $(e_n)$
 on a Hilbert space $X$ with corresponding eigenvalues $\la_n$.
If $Y=\CC$ and if
\[
   \mu = \sum_n |C e_n|^2 \delta_{-\la_n}
\]
is a Carleson measure on $\CC_+$, then $C$ is finite-time admissible
for the solution family $(S(t))_{t\ge 0}$.
\end{corollary}

A nice sufficient condition for admissibility for Cauchy
problems is given by Zwart \cite{Zwart:logplus}. For the convenience of
the reader we reproduce it here:

\begin{theorem}[Zwart] \label{thm:zwart}
Let $A$ be the infinitesimal generator of
an exponentially stable $C_0$-semigroup $T(t)_{t\ge 0}$ on the Hilbert
space $H$ and let $C:  X_1 \to Y$ be bounded, where $Y$ is another
Hilbert space. If for some  $\al>\einhalb$, 
\begin{equation} \label{eq:weiss-log-abschaetzung-zwart} 
\sup_{\Re \la>0}   \bignorm{ \,\ (1+\logplus \Re \la)^\al (\Re(\la))^\einhalb C
      R(\la, A) \,} < \infty, 
\end{equation}
then $C$ is a finite-time admissible observation operator.
\end{theorem}

Notice that the condition (\ref{eq:weiss-log-abschaetzung-zwart}) can be
reformulated by saying that in the sense of Evans, Opic and Pick (see
\cite{EvansOpicPick1996,EvansOpic,EvansOpicPick2002})
\[
\forall x\in X: \qquad \norm{ C R(\cdot, A)x }_{\einhalb,\infty. {\mathbb A}} < \infty
\]
where ${\mathbb A}=(0,\al)$, see also Cobos, Frenandez-Cabrera and
Triebel \cite{CobosFernandezTriebel} for logarithmic type
interpolation functors.

\begin{corollary}
Let in addition to the assumptions of Theorem~\ref{thm:zwart}, $a$ be a
$1$-regular and sectorial kernel of type $<\pihalbe$. Then $C$ is 
finite-time admissible for the solution family $\bigl(S(t)\bigr)_{t\ge0}$.
\end{corollary}

In some situations, the condition of sectoriality of angle $<\pihalbe$
in the above corollary may be inconvenient. Under merely $1$-regularity
one can also obtain admissibility by the following direct argument.

\begin{theorem}\label{thm:hans}
  Let $S(\cdot)$ be a bounded solution family to
  (\ref{eq:volterra-system}) with a $1$--regular kernel function $a$.
  Let $C: X_1 \to Y$ be bounded and assume that for some
  $\al>\einhalb$,
\begin{equation} \label{eq:weiss-log-abschaetzung} \sup_{r>0}
    \bignorm{ \,\ (1+\logplus r)^\al r^\einhalb C H(r) \,} < \infty.
\end{equation}
Then $C$ is finite-time admissible for $(S(t))_{t\ge 0}$.
\end{theorem}
Note that the exponent $\al>\einhalb$ is  optimal in the sense that
for $\al<\einhalb$ it is even wrong in the case $a\equiv 1$, see
\cite{JacobPartingtonPott:zero-class}. About the case $\al=\einhalb$
nothing is known at the moment.

\begin{proof}
  Let $\la\in\CC_+$ and let $\varphi$ such that $\la = |\la| e^{2
    i\varphi}$.  Then, by resolvent identity,
\begin{eqnarray*}
&&    (1+(\logplus( \Re \la))^\al \la^\einhalb C H(\la)  \\
&=&   (1+(\logplus(\Re\la))^\al \la^{-\einhalb} C \frac{1}{\ahut(\la)} R(\frac{1}{\ahut(\la)}, A)  \\
&=&   (1+\logplus |\la|)^\al |\la|^\einhalb C H(|\la|) \; e^{-i
\varphi } \frac{\ahut(|\la|)} {\ahut(\la)} \Bigl[ I + \bigl(
       \frac{1}{\ahut(|\la|)}-\frac{1}{\ahut(\la)}\bigr)
       R(\frac{1}{\ahut(\la)}, A) \Bigr]  \\
&=&    (1+\logplus |\la|)^\al |\la|^\einhalb C H(|\la|) \; e^{-i  \varphi } \Bigl[
       I + (1-\frac{\ahut(|\la|)} {\ahut(\la)}) A R(\frac{1}{\ahut(\la)}, A)\Bigr].
\end{eqnarray*}
By \cite[Lemma 8.1]{Pruess:Volterra}, 
$c^{-1} \leq \bigl|\ahut(|\la|) / \ahut(\la) \bigr| \leq c$ for some
$c>0$.  This yields uniform boundedness of expression in brackets and
so the assumed estimate (\ref{eq:weiss-log-abschaetzung}) gives
\begin{equation}  \label{eq:kleines-realteillemma}
\norm{ \la \mapsto C H(r{+}\la) }_{H^\infty(\CC_+)} \le K (1+\logplus r)^{-\al} r^{-\einhalb}.
\end{equation}
Since $(S(t))_{t\ge 0}$ is bounded,
\[
    \norm{\la \mapsto  H(r{+}\la)x }_{H^2(\CC_+)}
=   \norm{ e^{-rt} S(t) x }_{H^2(\CC_+)}
\le K r^{-\einhalb}\, \norm{x} \qquad \forall r>0
\]
and together with (\ref{eq:kleines-realteillemma}), we infer
\begin{equation} \label{eq:hans-eins}
    \norm{ \la \mapsto C H(r{+}\la)^2 x }_{H^2(\CC_+)}
\le \frac{K}{(1+\logplus r)^\al r}\, \norm{x}
\qquad \forall r>0.
\end{equation}
Moreover, the estimate
\[
     \bignorm{ \la\mapsto \frac1{r{+}\la} C H(r{+}\la) x}_{H^2(\CC_+)}
\le  \bignorm{ \la\mapsto C H(r{+}\la)x }_{H^\infty(\CC_+)}
     \bignorm{ \la\mapsto \frac1{r{+}\la} }_{H^2(\CC_+)}
\]
implies
\begin{equation}   \label{eq:hans-zwei}
    \bignorm{\la\mapsto \frac1{r{+}\la} C H(r{+}\la)x }_{H^2(\CC_+)}
\le \frac{K}{(1+\logplus r)^\al r} \norm{x} \qquad \forall r>0.
\end{equation}
Since $\frac{d}{d\la} H(\la) = \Bigl(\frac{ \la
  \ahut'(\la)}{\ahut(\la)}\Bigr) H(\la)^2 - \frac1{\la} \bigl(1+
\frac{ \la \ahut'(\la)}{\ahut(\la)}\bigr)H(\la)$ we infer from
(\ref{eq:hans-eins}) and (\ref{eq:hans-zwei}) that
\[
    \Bignorm{\mu \mapsto  \frac{d}{d\mu} C H(r+\mu) x}_{H^2(\CC_+)}
\le \frac{K}{(1+\logplus r)^\al r} \norm{x} \qquad \forall r>0.
\]
Finally, (inverse) Laplace transform yields
\[
     \Bignorm{  t\mapsto r t e^{-r t} C S(t)x }_{L^2(\RR_+)}
 \le \frac{K}{(1+\logplus r)^\al} \norm{x}\qquad \forall r>0
\]
and that is the estimate we need in the following dyadic decomposition
argument: notice that $x e^{-x} \ge 2 e^{-2}$ for $x\in [1,2]$. Fix
some $t_0>0$.  Then,
\begin{eqnarray*}
      \int_0^{t_0} \bignorm{ C S(t)x }^2 \,dt
& = & \sum_{n=1}^\infty \int_{t_0 2^{-n}}^{t_0 2^{-n+1}} \bignorm{ C S(t)x }^2 \,dt \\
&\le& \tfrac{e^4}4 \sum_{n=1}^\infty \int_{t_0 2^{-n}}^{t_0 2^{-n+1}}
        \bignorm{ t 2^n t_0^{-1} e^{t 2^n t_0^{-1}} C S(t)x }^2 \,dt \\
&\le& K \sum_{n=1}^\infty \frac1{ (1+\logplus(2^n t_0^{-1}))^{2\al} } \norm{x}^2 \le K \norm{x}^2.
\end{eqnarray*}
\end{proof}

\section{Example}

In this section we apply the obtained results to time-fractional
diffusion equation of distributed order.

Let $A$ generate an exponentially stable strongly continuous
semigroup $(T(t))_{t\ge 0}$. For $\omega >0$, $0<2\alpha <\beta\le
1$ we study a time-fractional diffusion equation of distributed
order of the form
\begin{eqnarray}
\omega D^\alpha_t x(t)
+ D^\beta_t x(t)&=& Ax (t),\quad t\ge 0,\label{eqnfractional}\\
x(0) &=& x_0,\nonumber
\end{eqnarray}
where $D^\alpha_t x = \bigl(-\frac\partial{\partial t}\bigr)^\alpha
x$ denotes the Caputo derivative of $x$, given by the Phillips
functional calculus of the right shift semigroup, that is,
\[ D^\gamma_t x(t)=\frac{1}{\Gamma(1-\gamma)} \int_0^t (t-s)^{-\gamma}x'(s)\,ds.\]
for $\gamma\in(0,1)$. Time-fractional diffusion equations of
distributed order have attracted attention as a possible tool for
the description of anomalous diffusion and relaxation phenomena in
many areas such as turbulence, disordered medium, intermittent
chaotic systems, mathematical finance and stochastic mechanics. For
further information on time-fractional diffusion equations of
distributed order we refer the reader to
\cite{Atanackovic,AtanackovicBudincevic,BayleyTorvik2,BayleyTorvik1,CarpinteriMainardi,ChechkinGorenflo,Jumarie,
MainardiPagnini}.

Using the Laplace transform equation (\ref{eqnfractional}) is
equivalent to
\[ x(t) = x_0 + \int_0^t a(t{-}s)\, A x(s)\, ds,\]
where \[ a(t)= t^{\beta-2\alpha-1} E_{\beta-\alpha, \beta-2\alpha}(
-\omega t^{\beta-\alpha}).\] Here $E_{\gamma,\delta}$, where
$\gamma, \delta>0$, denotes the Mittag-Leffler function
\[ E_{\gamma,\delta}(z) =\sum_{k=0}^\infty
\frac{z^k}{\Gamma(\gamma k+\delta)}\] The Laplace transformation of
the kernel $a$ is given by \[
\hat{a}(\lambda)=\frac{\lambda^{-\alpha}}{\omega
+\lambda^{\beta-\alpha}}.\] Thus the kernel $a$ satisfies the
assumption of Theorem \ref{thm:theo1}.

\medskip

We note that this example does e.g. not satisfy the assumption of
\cite[Theorem 3.10]{hjpp08} due to the 'mixed' growth conditions near
infinity and the origin, such that, even when $A$ is the Dirichlet
Laplacian on a bounded domain, the latter result cannot be used to
guarantee admissibility whereas a 'standard' Carlseon measure
reduction using Corollary~\ref{cor:carleson} still applies.

\section*{Acknowledgement}

This research was done at the Mathematisches Forschungsinstitut
Oberwolfach during a stay within the Research in Pairs Programme from
June 15 to June 28, 2008. We would like to thank the MFO for excellent
working conditions.

\def\SUBMITTED{Submitted}
\def\TOAPPEAR{To appear}
\def\PREPARATION{In preparation}

\def\cprime{$'$}

\end{document}